\documentclass[10pt,reqno]{article}
\usepackage{latexsym, amsmath, amssymb, amsthm, a4, epsfig}
\usepackage{graphicx}
\usepackage[usenames]{color}

\newtheorem{theorem}{Theorem}[section]
\newtheorem{cor}[theorem]{Corollary}
\newtheorem{lemma}[theorem]{Lemma}
\newtheorem{prop}[theorem]{Proposition}

\newcommand{\nm}{\noalign{\smallskip}}
\newcommand{\ds}{\displaystyle}
\newcommand{\p}{\partial}

\newcommand{\eqnref}[1]{(\ref {#1})}

\newcommand{\Rbb}{\mathbb{R}}

\newcommand{\la}{\langle}
\newcommand{\ra}{\rangle}

\newcommand{\Ccal}{\mathcal{C}}

\newcommand{\Hcal}{\mathcal{H}}
\newcommand{\Lcal}{\mathcal{L}}
\newcommand{\Kcal}{\mathcal{K}}

\newcommand{\Qcal}{\mathcal{Q}}
\newcommand{\Rcal}{\mathcal{R}}
\newcommand{\Scal}{\mathcal{S}}

\def\BF{{\bf F}}

\newcommand{\Ga}{\alpha}

\newcommand{\Gd}{\delta}
\newcommand{\Ge}{\epsilon}

\newcommand{\Gvf}{\varphi}
\newcommand{\Gg}{\gamma}

\newcommand{\Gl}{\lambda}

\newcommand{\Gs}{\sigma}

\newcommand{\Gj}{\tau}

\newcommand{\Go}{\omega}

\newcommand{\Gy}{\psi}

\newcommand{\GD}{\Delta}

\newcommand{\GG}{\Gamma}

\newcommand{\GO}{\Omega}

\newcommand{\beq}{\begin{equation}}
\newcommand{\eeq}{\end{equation}}

\def\ol{\overline}

\numberwithin{equation}{section}
\numberwithin{figure}{section}

\begin{document}
\title{Plasmon resonance with finite frequencies: a validation of the quasi-static approximation for diametrically small inclusions}

\author{Kazunori Ando\thanks{Department of Mathematics, Inha University, Incheon 402-751, S. Korea. Email: {\tt ando@inha.ac.kr, hbkang@inha.ac.kr}}, \and Hyeonbae Kang\footnotemark[1] \and Hongyu Liu\thanks{Department of Mathematics, Hong Kong Baptist University, Kowloon, Hong Kong SAR. Email: {\tt hongyuliu@hkbu.edu.hk}}}

\maketitle

\begin{abstract}
We study resonance for the Helmholz equation with a finite frequency in a plasmonic material of negative dielectric constant in two and three dimensions. We show that the quasi-static approximation is valid for diametrically small inclusions. In fact, we quantitatively prove that if the diameter of a inclusion is small compared to the loss parameter, then resonance occurs exactly at eigenvalues of the Neumann-Poincar\'e operator associated with the inclusion.
\end{abstract}

\noindent{\footnotesize {\bf AMS subject classifications}.}


\section{Introduction}\label{sec:intro}

Resonance phenomena is often observed in nanoscale particles whose material has a negative dielectric permittivity with a large wavelength in comparison with particle dimensions both experimentally and numerically \cite{PhysRevB.72.155412}.
It is known that such a resonance only occurs at certain frequencies.
Noble metals such as gold and silver show a negative permittivity \cite{palik1998handbook}, and are called plasmonic materials. Recently, there has been considerable interest in the plasmon resonance and its various applications including invisibility cloaking, biomedical imaging and medical therapy; see, e.g., \cite{Acm13,Ack13,ADM,AnKa14,MR3195185,JLSS,KLO,Klsap,LLL,Min06,PhysRevB.72.155412} and references therein.

It is known (see, e.g., \cite{AnKa14, MR3195185}) that in the quasi-static limit the plasmon resonance occurs at the eigenvalues of the Neumann-Poincar\'e operator associated with the inclusion. To be more precise, let $D$ be a bounded simply connected domain in $\Rbb^d$ ($d = 2, 3$) whose boundary $\p D$ is $\Ccal^{1, \Ga}$ for some $0<\alpha<1$.
Suppose that $D$ is occupied with a plasmonic material which has a dielectric constant $\Ge_c + i \Gd$, where $\Ge_c < 0$ and $\Gd > 0$ is the dissipation, and that the matrix $\Rbb \setminus \overline{\GO}$ has a dielectric constant $\Ge_m > 0$.
Hence, the distribution of the dielectric constant is given by
\beq
  \Ge_D =
  \begin{cases}
    \Ge_c + i \Gd, & \text{ in }\ D, \\
    \Ge_m, & \text{ in }\ \Rbb \setminus \ol{D}.
  \end{cases} \label{GeD}
\eeq
The dielectric equation in the quasi-static limit is given by
\beq
  \nabla \cdot \Ge_D \nabla u_\Gd = f.
\eeq
It is proved (e.g., \cite{AnKa14}) that when the source $f$ is given by the polarizable dipole $a \cdot \nabla \Gd_z$, the resonance occurs exactly when $\Gl(\Ge_c/\Ge_m)$ is an eigenvalue of the the Neumann-Poincar\'e (NP) operator associated with $D$ (see the next section for the definition and spectral properties of the NP operator), in other words, $\| \nabla u_\Gd \|_{L^2(D)} \to \infty$ as $\Gd \to \infty$. Here,
\beq
\Gl(t) := \frac{t+1}{2 ( t-1 )}.
\eeq
When $\Gl(\Ge_c/\Ge_m)$ is an eigenvalue of the the NP operator, $\Ge_c/\Ge_m$ is called the plasmon eigenvalue \cite{MR3195185}.

In this paper, we consider plasmon resonance for the Helmholtz operator $\nabla \cdot \Ge_D \nabla + \Go_0^2$, when $D$ is a diametrically small inclusion such as a nano-scale particle. Here $\Go_0$ represents the non-zero (but fixed) frequency, and the parameters $\Ge_c$ and $\Gd$ are determined by $\Go_0$. We show that if the diameter $s$ of $D$ is much smaller than the dissipation parameter $\Gd$, then the resonance occurs exactly when $\Gl(\Ge_c/\Ge_m)$ is an eigenvalue of the NP operator on $D$, like the quasi-static limit case. So the result of this paper can be regarded as a validation of the quasi-static approximation for small inclusions. It is worth mentioning that a different validation of quasi-static approximation is proved in \cite{ADM} by showing that the small volume asymptotic expansion of the far field for the Maxwell system holds away from the eigenvalues of the NP operator.

To describe results of this paper in a quantitative manner, let $D=s\GO$, and let after scaling
\beq
  \Ge_\GO =
  \begin{cases}
    \Ge_c + i \Gd, & \text{ in }\ \GO, \\
    \Ge_m, & \text{ in }\ \Rbb \setminus \overline{\GO}.
  \end{cases} \label{eq:dielectric_const2}
\eeq
We then consider
\beq
  \nabla \cdot \Ge_\GO \nabla u_\Gd + s^2 \Go_0^2 u_\Gd = a \cdot \nabla \Gd_{z} \quad \text{in }\ \Rbb^d \label{eq:Helmholtz2}
\eeq
satisfying the Sommerfeld radiation condition
\beq
  \left\vert \frac{\p u_\Gd}{\p r} - i \Go\varepsilon_m^{-1/2} u_\Gd \right\vert \le C r^{- \left( d + 1 \right) / 2}\ \text{ as }\ r = \left\vert x \right\vert \to \infty, \label{eq:radiation_condition}
\eeq
where $a \in \Rbb^d$ is a constant vector and $\Gd_z$ is the Dirac mass at $z \in \Rbb^d \setminus \overline{\GO}$. We characterize the resonance by the blow-up of $\left\Vert \nabla u_\Gd \right\Vert_{L^2(\GO)}$:
\beq
  \left\Vert \nabla u_\Gd \right\Vert_{L^2(\GO)} \to \infty\ \text{ as }\ \Gd \to +0. \label{eq:solution_blow-up}
\eeq
where $u_\Gd$ is the solution to \eqref{eq:Helmholtz2}.

We show that if $s$ is much smaller than $\Gd$, more precisely, if $s \Gd^{-1} \ll 1$ in three dimensions, and $s^2 |\ln s| \Gd^{-1} \ll 1$ in two dimensions, then \eqnref{eq:solution_blow-up} takes place if and only if $\Gl(\Ge_c/\Ge_m)$ is an eigenvalue of the NP operator on $\GO$. Moreover, if $\Gl(\Ge_c/\Ge_m)$ is an eigenvalue, we obtain a quantitative estimate
\beq\label{quanti}
\left\Vert \nabla u_\delta \right\Vert_{L^2(\Omega)} \approx \delta ^{- 1} \quad\text{as } \delta \to +0
\eeq
for most $z$ (the location of the dipole source). See Theorem \ref{mainthm3D} and Theorem \ref{mainthm2D} for precise statements. It is worth mentioning that the spectrums of the NP operators on $D$ and on $\GO$ are the same.

The rest of this paper is organized as follows. In section \ref{sec:prelim} we review spectral properties of the NP operator.
Section \ref{sec:asymp_exp} is to derive necessary asymptotic formula for the Helmholtz operator at low frequencies and estimates for the $H^1$-norm of the solution.  The main results in three and two dimensions are presented and proved in subsection \ref{sec:3d} and  subsection \ref{sec:2d}, respectively.

While writing this paper (after completing major work) we received the paper \cite{AMRZ15} from Habib Ammari. There an asymptotic formula for the solution similar to \eqnref{solGvf3D} is derived in three dimensions when there are multiple small inclusions, using the same method as in this paper (the spectral properties of the NP operator). Then the formula is used to study enhancement of scattering and absorption, and super-resonance by plasmonic particles. Here in this paper we use the asymptotic formula to show resonance quantified by \eqnref{quanti}.

\section{Preliminaries}\label{sec:prelim}

Let $\GO$ be a bounded domain with the Lipschitz boundary in $\Rbb^d$, $d=2,3$. Throughout this paper $H^s(\p\GO)$ denotes the $L^2$-Sobolev space on $\p\GO$ whose norm is expressed as $\left\Vert \cdot \right\Vert_s$.
We denote by $\left\la \cdot, \cdot \right\ra$ the duality pairing of $H^{- 1 / 2}(\p \GO)$ and $H^{1 / 2}(\p \GO)$.
Let $H_0^{- 1 / 2}(\p \GO)$ be the space of $\Gy \in H^{- 1 / 2}(\p \GO)$ satisfying $\left\la \Gy, 1 \right\ra = 0$.

Let $\GG(x)$ be the fundamental solution to the Laplacian on $\Rbb^d$, $d = 2, 3$:
\beq
  \GG(x) =
  \begin{cases}
    \ds \frac{1}{2 \pi} \ln {\left\vert x \right\vert}, & d = 2, \\
    - \ds \frac{1}{4 \pi} \left\vert x \right\vert^{- 1}, & d = 3.
  \end{cases} \nonumber
\eeq
The single layer potential of $\Gvf \in H^{- 1 / 2}(\p \GO)$ for the Laplacian is defined by
\beq
  \Scal [\Gvf] (x) = \int_{\p \GO} \GG(x - y)\Gvf(y) d\sigma(y), \quad x \in \Rbb^d. \nonumber
\eeq
It is well known (see e.g. \cite{MR2327884}) that the following jump formula holds:
\beq
\p_\nu \Scal [\Gvf] \big|_\pm (x) = \left( \pm 1/2 I + \Kcal^* \right) [\Gvf] (x), \quad x \in \p \GO, \label{eq:jump_formula}
\eeq
where $\Kcal^*$ is the Neumann-Poincar\'{e} (NP) operator defined by
\beq
  \Kcal^* [\Gvf] (x) = \mathrm{p.v.} \int_{\p \GO} \p_{\nu_x} \GG(x - y) \Gvf(y) d\sigma(y), \quad x \in \p \GO. \label{eq:NP operator}
\eeq
Here $\p_\nu$ denotes the outward normal derivative, the subscripts $\pm$ the limit (to $\p \GO$) from outside and inside of $\GO$, respectively, and p.v. the Cauchy principal value.

It is proved in \cite{MR2308861} (see also \cite{Ka15}) that the NP operator $\Kcal^*$ can be symmetrized using Plemelj's symmetrization principle:
\beq
  \Scal \Kcal^* = \Kcal \Scal. \label{eq:Plemelj}
\eeq
In fact, if we define a new inner product on $H_0^{- 1 / 2}(\p \GO)$ by
\beq
  \left\la \Gvf, \Gy \right\ra_{\Hcal^*} := - \left\la \Gvf, \Scal [\Gy] \right\ra, \label{eq:InnerProd_H0*}
\eeq
where the right hand side of \eqref{eq:InnerProd_H0*} is well-defined since $\Scal$ maps $H^{- 1 / 2}(\p \GO)$ into $H^{1 / 2}(\p \GO)$,
then $\Kcal^*$ is self-adjoint with respect to this inner product. Let $\Hcal_0^*$ be the space $H_0^{- 1 / 2}(\p \GO)$ equipped with the inner product $\left\la \cdot, \cdot \right\ra_{\Hcal^*}$ and $\Vert \cdot \Vert_{\Hcal^*}$ be the induced norm. It is known (see \cite{KKLS}) that $\Vert \cdot \Vert_{\Hcal^*}$ is equivalent to the norm $\left\Vert \cdot \right\Vert_{- 1 / 2}$:
\beq
  \left\Vert \Gvf \right\Vert_{\Hcal^*} \approx \left\Vert \Gvf \right\Vert_{- 1 / 2} \label{eq:NormEquiv_H*}
\eeq
for $\Gvf \in H_0^{- 1 / 2}(\p \GO)$. Here and throughout this paper $A \lesssim B$ means $A \le C B$ for some constant $C$ independent of parameters involved; $A \approx B$ means that $A \lesssim B$ and $A \gtrsim B$.

There is a nontrivial $\Gvf_0 \in H^{-1/2}(\p \GO)$ such that
\beq\label{Gvf0}
\Kcal^* [\Gvf_0] = \frac{1}{2} \Gvf_0,
\eeq
We note that $\Scal [\Gvf_0]$ is constant, say $c_0$, in $\GO$. In three dimensions, $c_0 \neq 0$, and hence $\Scal: H^{-1/2}(\p \GO) \to H^{1/2}(\p \GO)$ is invertible. However, there are domains $\GO$ in two dimensions such that $c_0=0$ (see \cite{MR769382}), which means $\Scal$ is not invertible in general.
We introduce a variance of the single layer potential, denoted by $\tilde{\Scal}$, by $\tilde{\Scal}=\Scal$ if $c_0 \neq 0$ , and if $c_0 =0$, then
\beq
  \tilde{\Scal} [\Gvf] =
  \begin{cases}
    \Scal [\Gvf], &\ \text{ if }\ \left\la \Gvf, 1 \right\ra = 0, \\
    1, &\ \text{ if }\ \Gvf = \Gvf_0.
  \end{cases} \nonumber
\eeq
Then $\tilde{\Scal}$ is a bijection from $H^{-1/2}(\p \GO)$ to $H^{1 / 2}(\p \GO)$.
Moreover, we have an extension of \eqref{eq:Plemelj}:
\beq
  \tilde{\Scal} \Kcal^* = \Kcal \tilde{\Scal} \label{eq:ex_Plemelj}
\eeq
which enables us to extend the inner product \eqref{eq:InnerProd_H0*} to $H^{- 1 / 2}(\p \GO)$:
\beq
  \left\la \Gvf, \Gy \right\ra_{\Hcal^*} := - \la \Gvf, \tilde{\Scal}_{\p \GO} [\Gy] \ra. \label{eq:InnerProd_H*}
\eeq
We denote by $\Hcal^*$ the space $H^{- 1 / 2}(\p \GO)$ equipped with the inner product \eqref{eq:InnerProd_H*}. Then
the symmetrization principle \eqref{eq:ex_Plemelj} makes $\Kcal^*$ self-adjoint on $\Hcal^*$.
We emphasize that the norm equivalence \eqref{eq:NormEquiv_H*} is valid for $\Gvf \in H^{- 1 / 2}(\p \GO)$.

The spectrum $\sigma(\Kcal^*)$ of the NP operator lies in $( - 1 / 2, 1 / 2 ]$. Moreover, $\Kcal_{\p \GO}^*$ is a compact operator on $\Hcal^*$, when $\p \GO$ is $\Ccal^{1, \Ga}$ for some $\alpha>0$.
Therefore, we have the spectral decomposition of $\Kcal^*$ on $\Hcal^*$:
\beq
  \Kcal^* = \frac{1}{2} \Gvf_0 \otimes \Gvf_0 + \sum_{n = 1}^\infty \Gl_n \Gvf_n \otimes \Gvf_n = \sum_{n = 0}^\infty \Gl_n \Gvf_n \otimes \Gvf_n, \label{eq:SpecResol_NP*}
\eeq
where $\Gvf_n \in \Hcal^*$ is an eigenvector of $\Kcal^*$ corresponding to the eigenvalue $\Gl_n \in \Rbb$ (counting multiplicities), with $ 1 / 2= \Gl_0 > \left\vert \Gl_1 \right\vert \ge \left\vert \Gl_2 \right\vert \cdots \ge \left\vert \Gl_n \right\vert \ge \cdots \to 0$ as $n \to \infty$.
We note that $\{ \Gvf_n \}_{n = 0}^\infty$ is chosen to be an orthonormal basis on $\Hcal^*$ ($\Gvf_0$ is normalized so that $\| \Gvf_0 \|_{\Hcal^*}=1$).

Define an inner product
\beq
  \left\la f, g \right\ra_{\Hcal} := - \la f, \tilde{\Scal}^{- 1} [g] \ra
\eeq
on $H^{1/2}(\p \GO)$, and denote by $\Hcal$ the space $H^{1 / 2}(\p \GO)$ equipped with the inner product $\la \cdot, \cdot \ra_{\Hcal}$.
Then $\tilde{\Scal}$ is a unitary operator from $\Hcal^*$ to $\Hcal$, and hence $\{ \tilde{\Scal}[\Gvf_n], \ n=0,1,\ldots \}$ is an orthonormal basis of $\Hcal$. Let $\Hcal_0$ be the subspace of $\Hcal$ spanned by $\{ \tilde{\Scal}[\Gvf_n], \ n=1,\ldots \}$. Then $\Scal: \Hcal^*_0 \to \Hcal_0$ is a bijection. We emphasize that the norm $\left\Vert \cdot \right\Vert_{\Hcal}$ is equivalent to $\left\Vert \cdot \right\Vert_{1 / 2}$.

For $\Gvf \in \Hcal^*$, we write
\beq\label{Gvfhat}
\hat{\Gvf}(n):= \la \Gvf, \Gvf_n \ra_{\Hcal^*}, \quad n=0,1,2, \ldots,
\eeq
so that
\beq\label{specdec1}
\Gvf= \sum_{n=0}^\infty \hat{\Gvf}(n) \Gvf_n (= \hat{\Gvf}(0) \Gvf_0 + \Gvf^\prime), \quad \| \Gvf \|_{\Hcal^*}^2 = |\hat{\Gvf}(0)|^2 + \| \Gvf^\prime \|_{\Hcal^*}^2 .
\eeq
For $f \in \Hcal$, we define
\beq
\check{f}(n):= \la f, \tilde{\Scal}[\Gvf_n] \ra_{\Hcal}, \quad n=0,1,2, \ldots,
\eeq
so that
\beq\label{specdec2}
f= \sum_{n=0}^\infty \check{f}(n) \tilde{\Scal}[\Gvf_n] (= \check{f}(0) \tilde{\Scal}[\Gvf_0]+ f^\prime), \quad \| f \|_{\Hcal}^2 = |\check{f}(0)|^2 + \| f^\prime \|_{\Hcal}^2 .
\eeq
We refer to \cite{AnKa14} and references therein for more details on the preliminaries presented in this section.

Finally, we denote by $\Lcal(X, Y)$ the space of bounded linear operators from a Banach space $X$ to a Banach space $Y$; in particular, $\Lcal(X)$ is the space of  bounded linear operators on a Banach space $X$.

\section{Asymptotic expansion at low frequencies} \label{sec:asymp_exp}

Let $\Go=s\Go_0$ from now on to make notation short.
A fundamental solution $\GG^\Go(x)$ to the Helmholtz operator $\GD + \Go^2$ in $\Rbb^d$ is a solution of
\beq
  (\GD + \Go^2) \GG^\Go = \Gd_0, \label{eq:HelmholtzOp}
\eeq
where $\Gd_0$ is the Dirac function at $0$.
Among solutions to \eqref{eq:HelmholtzOp}, we seek a solution satisfying the Sommerfeld radiation condition
\beq
  \left\vert \frac{\p \GG^\Go}{\p r} - i \Go \GG^\Go \right\vert \le C r^{- (d + 1) / 2}\ \text{ as }\ r = \left\vert x \right\vert \to \infty. \label{eq:rad_cond}
\eeq
Then, it is given by
\beq
  \GG^\Go(x) =
  \begin{cases}
    - \ds \frac{i}{4} H_0^1(\Go \left\vert x \right\vert) \quad & \text{if } d = 2, \\
    \nm
    - \ds \frac{1}{4 \pi} \frac{e^{i \Go \left\vert x \right\vert}}{\left\vert x \right\vert} & \text{if } d = 3,
  \end{cases} \label{eq:FundSol_Helmholtz}
\eeq
where $H_0^1(z)$ is the Hankel function of the first kind of order $0$.

For the subsequent use, we consider the asymptotic expansion of the fundamental solution $\GG^\Go(x)$ as $\Go \to +0$. When $d = 2$, we recall the behavior of the Hankel function $H_0^1(z)$ near $z = 0$ (see, e.g., \cite{Leb72}):
\beq
  - \frac{i}{4} H_0^1(\Go \left\vert x \right\vert) = \frac{1}{2 \pi} \ln {\left\vert x \right\vert} + \Gj + \sum_{n = 1}^\infty \left( b_n \ln \left( \Go \left\vert x \right\vert \right) + c_n \right) \left( \Go \left\vert x \right\vert \right)^{2 n}, \label{eq:FundSol_exp2}
\eeq
where
\beq
  b_n = \frac{\left( - 1 \right)^n}{2 \pi} \frac{1}{2^{2 n} \left( n! \right)^2}, \quad
  c_n = - b_n \left( \Gg - \ln {2} - \frac{\pi i}{2} - \sum_{j = 1}^n \frac{1}{j} \right) \nonumber
\eeq
and
\beq
  \Gj = \frac{1}{2\pi} \left( \ln {\Go} + \Gg - \ln {2} \right) - \frac{i}{4} \label{eq:tau}
\eeq
($\Gg$ is the Euler constant). So we have
\beq
  \GG^\Go(x) = \GG(x) + \Gj + \Go^2 \ln {\Go} K_2^\Go(x) \label{eq:10}
\eeq
as $\Go \to +0$ (see also \cite{MR2488135}). The definition of $K_2^\Go(x)$ is obvious.
When $d = 3$, one can easily see that
\beq
  - \frac{1}{4 \pi} \frac{e^{i \Go \left\vert x \right\vert}}{\left\vert x\right\vert} = - \frac{1}{4 \pi} \frac{1}{\left\vert x \right\vert} - \frac{i \Go}{4 \pi} \sum_{n = 1}^\infty \frac{\left( i \Go \left\vert x \right\vert \right)^{n - 1}}{n!}, \label{eq:FundSol_exp3}
\eeq
which implies that
\beq
  \GG^\Go(x) =  \GG(x) + \Go K_3^\Go(x). \label{eq:11}
\eeq

Let us observe a regularity property of the function $K_d^\Go(x)$ ($d=2,3$) for later purpose. Let $\Go_1$ be a small positive number. Then there is a constant $C$ independent of $\Go \le \Go_1$ such that
\beq\label{Kdest1}
\int_{\GO} \int_{\p\GO} |\p_x^\Ga K_d^\Go (x-y)|^2 d\Gs(y) \, dx \le C
\eeq
for all $\Ga=(\Ga_1, \ldots, \Ga_d)$ such that $|\Ga| \le 2$. Here $\p_x^\Ga$ is the partial derivative with respect to $x$.
Moreover, $\nabla K_3^\Go(x)$ gains $\Go$ and it holds that
\beq\label{Kdest2}
\frac{1}{\Go} \int_{\GO} \int_{\p\GO} |\p_x^\Ga \nabla_x K_3^\Go (x-y)|^2 d\Gs(y) \, dx \le C
\eeq
for all $|\Ga| \le 1$.

The single layer potential of $\Gvf \in H^{- 1 / 2}(\p \GO)$ for the Helmholz operator $\GD + \Go^2$ is defined by
\beq
  \Scal^\Go [\Gvf] (x) = \int_{\p \GO} \GG^\Go(x - y) \Gvf(y) d\sigma(y), \quad x \in \Rbb^d.
\eeq
We note that $\Scal^\Go [\Gvf] (x)$ satisfies the Sommerfeld radiation condition \eqref{eq:rad_cond} (see \cite{MR2488135}). Let $\Rcal_d^\Go$ ($d=2,3$) be the integral operator defined by $K_d^\Go$, namely,
\beq
  \Rcal_d^\Go [\Gvf] (x) = \int_{\p \GO} K_d^\Go(x - y) \Gvf(y) d\sigma(y), \quad x \in \Rbb^d.
\eeq
Then, we obtain from \eqnref{eq:10} and \eqnref{eq:11} that
\beq\label{Stwo}
\Scal^\Go =
\begin{cases}
\Scal + \Gj \left\la \cdot, 1 \right\ra + \Go^2 \ln\Go \Rcal_2^\Go \quad & \mbox{if } d=2, \\
\Scal + \Go \Rcal_3^\Go \quad & \mbox{if } d=3.
\end{cases}
\eeq

Analogously to \eqref{eq:jump_formula}, the following jump formula holds:
\beq
\p_\nu \Scal^\Go [\Gvf] \big|_\pm (x) = \left( \pm 1/2 I + \left( \Kcal^\Go \right)^* \right) [\Gvf] (x), \quad  x \in \p \GO, \label{eq:Helmholtz_jump_formula}
\eeq
where $\left( \Kcal^\Go \right)^*$ is defined by
\beq
  \left( \Kcal^\Go \right)^* [\Gvf] (x) = \int_{\p \GO} \p_{\nu_x} \GG^\Go(x - y) \Gvf(y) d\sigma(y), \quad x \in \p \GO. \nonumber
\eeq
For $d=2,3$, let
\beq
\Qcal_d^\Go[\Gvf](x):=
\begin{cases}
\p_\nu \Rcal_2^\Go[\Gvf](x), \quad d=2, \\
\frac{1}{\Go} \p_\nu \Rcal_3^\Go[\Gvf](x), \quad d=3,
\end{cases}
x\in \p\GO.
\eeq
Then, we have
\beq\label{Qtwo}
( \Kcal^\Go )^* =
\begin{cases}
\Kcal^* + \Go^2 \ln\Go \Qcal_2^\Go \quad & \mbox{if } d=2, \\
\Kcal^* + \Go^2 \Qcal_3^\Go \quad & \mbox{if } d=3.
\end{cases}
\eeq

We now investigate the mapping property of $\Rcal_d^\Go$ and $\Qcal_d^\Go$. By Cauchy-Schwartz inequality we see from \eqnref{Kdest1} that
$$
\| \Rcal_d^\Go[\Gvf] \|_{W^{2,2}(\GO)} \le C \| \Gvf \|_{L^2(\p\GO)}.
$$
We also see from \eqnref{Kdest2} that
$$
\frac{1}{\Go} \| \nabla \Rcal_3^\Go[\Gvf] \|_{W^{1,2}(\GO)} \le C \| \Gvf \|_{L^2(\p\GO)}.
$$
By trace theorem, $\Rcal_d^\Go$ maps $L^2(\p\GO)$ into $H^{3/2}(\p\GO)$, and $\Qcal_d^\Go$ maps $L^2(\p\GO)$ into $H^{1/2}(\p\GO)$.
By duality, $\Rcal_d^\Go$ maps  $H^{-3/2}(\p \GO)$ into $L^2(\p \GO)$, and $H^{-1/2}(\p \GO)$ into $H^{1}(\p \GO)$ by interpolation. Likewise we see that $\Qcal_d^\Go$ maps  $H^{-1/2}(\p \GO)$ into $L^2(\p \GO)$. We summarize these properties in the following lemma.

\begin{lemma} \label{regularity}
For a given small positive number $\Go_1$, there exists a constant $C$ independent of $\Go \le \Go_1$ such that
\beq\label{Rdreg}
\| \Rcal_d^\Go[\Gvf] \|_{1} \le C \| \Gvf \|_{-1/2}
\eeq
and
\beq\label{Qdreg}
\| \Qcal_d^\Go[\Gvf] \|_{0} \le C \| \Gvf \|_{-1/2}
\eeq
for all $\Gvf \in H^{-1/2}(\p\GO)$.
\end{lemma}

\begin{prop} \label{prop:3}
Let $\Gvf \in \Hcal^*$ and $\Gvf = \Gvf^\prime + \hat{\Gvf}(0) \Gvf_0$ be its orthogonal decomposition
where $\Gvf^\prime \in \Hcal_0^*$. The following estimates hold:
\begin{enumerate}
\item[(i)] If $d = 2$, then
\beq\label{eq:energy_estimate_2d}
\Vert \Gvf^\prime \Vert_{\Hcal^*}^2 - |\Go\ln\Go|^2 |\hat{\Gvf}(0)|^2 \lesssim \left\Vert \nabla \Scal_{\p \GO}^\Go [\Gvf] \right\Vert_{L^2(\GO)}^2 \lesssim \Vert \Gvf^\prime \Vert_{\Hcal^*}^2 +  |\Go\ln\Go|^2 |\hat{\Gvf}(0)|^2.
\eeq
\item[(ii)] If $d = 3$, then
\beq\label{eq:energy_estimate_3d}
\Vert \Gvf^\prime \Vert_{\Hcal^*}^2 - |\Go| |\hat{\Gvf}(0)|^2 \lesssim \left\Vert \nabla \Scal_{\p \GO}^\Go [\Gvf] \right\Vert_{L^2(\GO)}^2 \lesssim \Vert \Gvf^\prime \Vert_{\Hcal^*}^2 +  |\Go| |\hat{\Gvf}(0)|^2.
\eeq
\end{enumerate}
\end{prop}

\begin{proof}
We only prove \eqnref{eq:energy_estimate_2d} since three dimensional case can be proved in a similar way.

We have from Gauss's divergence theorem
\beq\label{Gauss}
\int_\GO \left\vert \nabla \Scal^\Go [\Gvf] \right\vert^2 dx + \int_\GO \Scal^\Go [\Gvf]\, \ol{\GD \Scal^\Go [\Gvf]} dx = \int_{\p \GO} \Scal^\Go [\Gvf] \, \ol{\p_\nu \Scal^\Go [\Gvf]|_- } \, d\Gs.
\eeq
Since $\GD \Scal^\Go [\Gvf]=- \Go^2 \Scal^\Go [\Gvf]$, we have
\beq
\int_\GO \left\vert \nabla \Scal^\Go [\Gvf] \right\vert^2 dx = \Go^2 \int_\GO \left| \Scal^\Go [\Gvf] \right|^2 dx + \int_{\p \GO} \Scal^\Go [\Gvf] \, \ol{\p_\nu \Scal^\Go [\Gvf]|_- } \, d\Gs. \label{eq:gauss_divergence}
\eeq
One can see from \eqnref{Stwo} and Lemma \ref{regularity} that
\beq\label{triest}
\int_\GO \left| \Scal^\Go [\Gvf] \right|^2 dx \lesssim |\ln \Go| \, \|\Gvf \|_{-1/2} \lesssim |\ln \Go| \, \|\Gvf \|_{\Hcal^*},
\eeq
since $|\Gj| \lesssim |\ln \Go|$. The last inequality holds because of \eqnref{eq:NormEquiv_H*}.

Using the jump formula \eqref{eq:Helmholtz_jump_formula} we have
$$
\int_{\p \GO} \Scal^\Go [\Gvf] \, \ol{\p_\nu \Scal^\Go [\Gvf]|_- } \, d\Gs
=  \int_{\p \GO} \Scal^\Go [\Gvf] \, \ol{\left( - 1/2I + \left( \Kcal^\Go \right)^* \right) [\Gvf]} \, d\Gs .
$$
One then see from \eqnref{Stwo} and \eqnref{Qtwo} that
\begin{align}
& \int_{\p \GO} \Scal^\Go [\Gvf] \, \ol{\p_\nu \Scal^\Go [\Gvf]|_- } \, d\Gs \nonumber \\
& = \int_{\p \GO} \Scal [\Gvf] \, \ol{\left( - 1/2I + \Kcal^* \right) [\Gvf]} \, d\Gs
 + \Gj \la \Gvf, 1 \ra \int_{\p \GO} \ol{\left( - 1/2I + \Kcal^* \right) [\Gvf]} \, d\Gs + \Go^2 \ln \Go E \label{formula1}
\end{align}
where
$$
E= \int_{\p \GO} \Rcal_2^\Go [\Gvf] \, \ol{\left( - 1/2I + \left( \Kcal^\Go \right)^* \right) [\Gvf]} \, d\Gs + \int_{\p \GO} \Scal^\Go [\Gvf] \, \ol{\Qcal_2^\Go [\Gvf]} \, d\Gs + \Gj \la \Gvf, 1 \ra \int_{\p \GO} \ol{\Qcal_2^\Go [\Gvf]} \, d\Gs.
$$

Using \eqnref{Stwo} and Lemma \ref{regularity} one can show that
\beq\label{Eest}
|E| \le C |\ln \Go| \| \Gvf \|_{\Hcal^*}^2
\eeq
for some constant $C$ independent of $\Go \le \Go_1$. In fact, we have from \eqnref{eq:tau}
\begin{align*}
|E| &\le \| \Rcal_2^\Go [\Gvf] \|_{1/2} \| ( - 1/2I + ( \Kcal^\Go )^* ) [\Gvf] \|_{-1/2} + \| \Scal^\Go [\Gvf] \|_{1/2} \| \Qcal_2^\Go [\Gvf] \|_{-1/2} + \tau \| \Gvf \|_{-1/2} \| \Qcal_2^\Go [\Gvf] \|_{0} \\
& \le C |\ln \Go| \| \Gvf \|_{-1/2}^2 .
\end{align*}

Since $\Kcal[1]=1/2$, we have
\beq\label{Kone}
\int_{\p \GO} \ol{\left( - 1/2I + \Kcal^* \right) [\Gvf]} \, d\Gs = \int_{\p \GO} \left( - 1/2I + \Kcal \right) [1] \, \ol{\Gvf} \, d\Gs =0.
\eeq
On the other hand, since $\Kcal^*[\Gvf_0]=1/2 \Gvf_0$, we have
$$
\int_{\p \GO} \Scal [\Gvf] \, \ol{\left( - 1/2I + \Kcal^* \right) [\Gvf]} \, d\Gs
= \int_{\p \GO} \Scal [\Gvf^\prime] \, \ol{\left( - 1/2I + \Kcal^* \right) [\Gvf^\prime]} \, d\Gs.
$$
Using $\Gvf^\prime = \sum_{n = 1}^\infty \hat{\Gvf}(n) \Gvf_n$, we have
$$
\int_{\p \GO} \Scal [\Gvf] \, \ol{\left( - 1/2I + \Kcal^* \right) [\Gvf]} \, d\Gs
= \sum_{n, m = 1}^\infty \left(- 1/2 + \Gl_l \right) \hat{\Gvf}(n) \ol{\hat{\Gvf}(m)} \int_{\p \GO} \Scal [\Gvf_j]  \, \ol{\Gvf_l} \, d\Gs .
$$
Since $\int_{\p \GO} \Scal [\Gvf_n]  \, \ol{\Gvf_m} \, d\Gs= -\la \Gvf_n, \Gvf_m \ra_{\Hcal^*}= -\Gd_{nm}$ (the Kronecker's delta), we have
\begin{align*}
\int_{\p \GO} \Scal [\Gvf] \, \ol{\left( - 1/2I + \Kcal^* \right) [\Gvf]} \, d\Gs
= \sum_{n=1}^\infty \left( \Gl_n - 1/2 \right) \left\vert \hat{\Gvf}(n) \right\vert^2 .
\end{align*}
So we have
\beq\label{formula2}
\left| \int_{\p \GO} \Scal [\Gvf] \, \ol{\left( - 1/2I + \Kcal^* \right) [\Gvf]} \, d\Gs \right| \approx \|\Gvf^\prime \|_{\Hcal^*}^2.
\eeq

Combining \eqnref{formula1}-\eqnref{formula2} we obtain
$$
\|\Gvf^\prime \|_{\Hcal^*}^2 - |\Go \ln \Go|^2 \|\Gvf \|_{\Hcal^*}^2 \lesssim \left| \int_{\p \GO} \Scal^\Go [\Gvf] \, \ol{\p_\nu \Scal^\Go [\Gvf]|_- } \, d\Gs \right| \lesssim \|\Gvf^\prime \|_{\Hcal^*}^2 + |\Go \ln \Go|^2 \|\Gvf \|_{\Hcal^*}^2 ,
$$
which together with \eqnref{Gauss} and \eqnref{triest} yields \eqnref{eq:energy_estimate_2d}.
\end{proof}

\section{Analysis of resonance}

From now on, we assume that $\Ge_m=1$ without loss of generality.

Set $k_m = \Go (=s\Go_0)$ and
\beq
  k_c^2 = \frac{\Go^2}{\Ge_c + i \Gd} , \ \ \Re {k_c} > 0, \ \ \Im {k_c} < 0. \nonumber
\eeq
Since
$$
  k_c = \Go \left( \Ge_c + i \Gd \right)^{- 1 / 2} \simeq - i \frac{\Go}{\sqrt{\left\vert \Ge_c \right\vert}} \left( 1 - i \frac{\Gd}{2\Ge_c} \right),
$$
we assume for simplicity
\beq
  k_c = - i \frac{\Go}{\sqrt{\left\vert \Ge_c \right\vert}} \left( 1 - i \frac{\Gd}{2 \Ge_c} \right). \label{eq:kcd}
\eeq
Then the problem \eqref{eq:Helmholtz2} can be written as
\beq\label{eq:transmission_3d}
  \begin{cases}
    \GD u_\Gd + k_c^2 u_\Gd = 0 \quad & \text{in }\ \GO, \\
    \GD u_\Gd + \Go^2 u_\Gd = a \cdot \nabla \Gd_{z} \quad & \text{in }\ \Rbb^d \setminus \ol{\GO}, \\
    \left. u_\Gd \right|_- - \left. u_\Gd \right|_+ = 0 & \text{on }\ \p \GO, \\
    \left( \Ge_c + i \Gd \right) \p_\nu u_\Gd |_- - \p_\nu u_\Gd |_+ = 0 & \text{on }\ \p \GO,
  \end{cases}
\eeq
under the Sommerfeld radiation condition \eqref{eq:radiation_condition}.

Let
\beq
F_z(x) := - a \cdot \nabla_x \GG^{\Go}(x - z).
\eeq
Then, the solution $u_\Gd$ can be represented as
\beq\label{eq:21}
  u_\Gd(x) =
  \begin{cases}
    \Scal^{k_c} [\Gvf_\Gd] (x), \quad & x \in \GO, \\
    F_z(x) + \Scal^{\Go} [\Gy_\Gd] (x), &x \in \Rbb^d \setminus \GO
  \end{cases}
\eeq
for some $\Gvf_\Gd, \psi_\Gd \in \Hcal^*$. In view of transmission conditions on $\p\GO$ (the third and fourth conditions in \eqnref{eq:transmission_3d}), $(\Gvf_\Gd, \psi_\Gd)$ should solve the following system of integral equations:
\beq\label{eq:system_integral_omega}
  \begin{cases}
    \Scal^{k_c} [\Gvf_\Gd] - \Scal^{\Go} [\Gy_\Gd] = F_z , \\
    ( \Ge_c + i \Gd ) \p_\nu \Scal^{k_c} [\Gvf_\Gd]|_-  - \p_\nu \Scal^{\Go} [\Gy_\Gd] |_+ = \p_\nu F_z,
  \end{cases}
  \text{ on } \p \GO.
\eeq

Let $X := \Hcal^* \times \Hcal^*$ and $Y := \Hcal \times \Hcal^*$, and define an operator $A_\Gd^s: X \to Y$ by
\beq
  A_\Gd^s =
  \begin{bmatrix}
    \Scal^{k_c} & - \Scal^{\Go} \\
   ( \Ge_c + i \Gd ) \p_\nu \Scal^{k_c}|_- & - \p_\nu \Scal^{\Go} |_+
  \end{bmatrix}. \label{eq:A_delta^omega}
\eeq
Then we can rewrite \eqref{eq:system_integral_omega} as
\beq\label{short}
  A_\Gd^s
  \begin{bmatrix}
    \Gvf_\Gd \\
    \Gy_\Gd
  \end{bmatrix}
  =
  \begin{bmatrix}
    F_z \\ \p_\nu F_z
  \end{bmatrix}.
\eeq
The solvability of \eqref{eq:system_integral_omega} is equivalent to the invertibility of $A_\Gd^s$.
We will investigate the behavior of the norm $\left( A_\Gd^s \right)^{- 1}$ as $\Gd \to +0$.

\subsection{Three dimensions}\label{sec:3d}

We deal with the three dimensional case first since it is easier.

We split $A_\Gd^s$ into two parts: $A_\Gd^s = A_\Gd + T_\Gd^s$, where
\beq
  A_\Gd =
  \begin{bmatrix}
    \Scal & - \Scal \\
    ( \Ge_c + i \Gd ) (-1/2I + \Kcal^*) & - (1/2I + \Kcal^*)
  \end{bmatrix}. \label{A_delta^0}
\eeq
Then we can infer from \eqnref{Stwo}, \eqnref{Qtwo} and Lemma \ref{regularity} that
\beq\label{TGd}
  \left\Vert T_\Gd^s \right\Vert_{\Lcal(X, Y)} \lesssim \Go.
\eeq

\begin{lemma}\label{thm3D}
For $f \in \Hcal$ and $g \in \Hcal^*$, the solution to
\beq\label{fg}
A_\Gd \begin{bmatrix} \Gvf \\ \Gy \end{bmatrix} = \begin{bmatrix} f \\ g \end{bmatrix}
\eeq
is given by
\beq\label{solGvf3D}
\Gvf= \sum_{n=0}^\infty \frac{\hat{g}(n) - (1/2+\Gl_n) \check{f}(n)}{( \Ge_c -1 ) (\Gl_n-\Gl(\Ge_c)) + i\Gd (\Gl_n-\frac{1}{2})} \Gvf_n
\eeq
and
\beq\label{solGy3D}
\Gy = \Gvf - \Scal^{-1}[f].
\eeq
\end{lemma}

\begin{proof}
The equation \eqnref{fg} can be written as
$$
\begin{cases}
    \Scal [\Gvf] - \Scal [\Gy] = f , \\
    ( \Ge_c + i \Gd ) (-1/2I + \Kcal^*)[\Gvf]  - (1/2I + \Kcal^*) [\Gy]  = g,
  \end{cases}
  \text{ on } \p \GO.
$$
Since $\Scal: \Hcal^* \to \Hcal$ is invertible in three dimensions, we have
\beq
\Gy = \Gvf - \Scal^{-1}[f].
\eeq
Substituting this into the second equation, we obtain
$$
\big( -1/2( \Ge_c + i \Gd +1 ) I + ( \Ge_c + i \Gd -1 )\Kcal^* \big) [\Gvf]  = g - (1/2I + \Kcal^*) \Scal^{-1}[f].
$$
We then use the spectral decomposition \eqnref{specdec1} to obtain
$$
\Gvf= \sum_{n=0}^\infty \frac{a_n}{-1/2( \Ge_c + i \Gd +1 ) + ( \Ge_c + i \Gd -1 ) \Gl_n} \Gvf_n
$$
where
$$
a_n = \hat{g}(n) - \la (1/2I + \Kcal^*) \Scal^{-1}[f], \Gvf_n \ra_{\Hcal^*}.
$$
Since $f= \sum_{j=0}^\infty \check{f}(j) \Scal [\Gvf_j]$, we have
\begin{align*}
\la (1/2I + \Kcal^*) \Scal^{-1}[f], \Gvf_n \ra_{\Hcal^*} = \sum_{j=0}^\infty \check{f}(j)  \la (1/2I + \Kcal^*)[\Gvf_j], \Gvf_n \ra_{\Hcal^*}
= (1/2+\Gl_n) \check{f}(n)  .
\end{align*}
This completes the proof.
\end{proof}

As a consequence of Theorem \ref{thm3D} we obtain the following corollary.

\begin{cor}\label{cor3D}
Suppose that $\Ge_c \neq -1$, and let $(\Gvf, \Gy)$ be the solution of \eqnref{fg}. Then the following hold for sufficiently small $\Gd$:
\begin{itemize}
\item[(i)] $\ds \Vert ( A_\Gd^0 )^{- 1} \Vert_{\Lcal(Y, X)} \lesssim \Gd^{-1}$.
\item[(ii)] If $\Gl(\Ge_c) \neq \Gl_n$ for any $n$, then $\ds \Vert ( A_\Gd^0 )^{- 1} \Vert_{\Lcal(Y, X)} \le C$ for some $C$ depending on $\Ge_c$.
\item[(iii)] If $\Gl(\Ge_c) = \Gl_n$ for some $n \neq 0$, then $\| \Gvf^\prime \|_{\Hcal^*} \gtrsim |a_n| \Gd^{-1}$, where $a_n= \hat{g}(n) - (1/2+\Gl_n) \check{f}(n)$.
\end{itemize}
\end{cor}

\begin{proof}
Since
$$
\frac{1}{|( \Ge_c -1 ) (\Gl_n-\Gl(\Ge_c)) + i\Gd (\Gl_n-\frac{1}{2})|} \lesssim \Gd^{-1},
$$
we have from \eqnref{solGvf3D} that
$$
\| \Gvf \|_{\Hcal^*}^2 \lesssim \Gd^{-2} \sum_{n=0}^\infty |\hat{g}(n) - (1/2+\Gl_n) \check{f}(n)|^2 \lesssim \Gd^{-2} (\| f \|_{\Hcal}^2 + \| g \|_{\Hcal^*}^2).
$$
We have from \eqnref{solGy3D} that
$$
\| \Gy \|_{\Hcal^*}^2 \lesssim \Gd^{-2} (\| f \|_{\Hcal}^2 + \| g \|_{\Hcal^*}^2).
$$
This proves (i).

Since $\Ge_c \neq -1$, $\Gl(\Ge_c) \neq 0$. If $\Gl(\Ge_c) \neq \Gl_n$ for any $n$, then $|\Gl_n-\Gl(\Ge_c)| \ge C$ for some $C>0$. So we have
$$
\frac{1}{|( \Ge_c -1 ) (\Gl_n-\Gl(\Ge_c)) + i\Gd (\Gl_n-\frac{1}{2})|} \lesssim 1,
$$
and hence
$$
\| \Gvf \|_{\Hcal^*}^2 \lesssim \sum_{n=0}^\infty |\hat{g}(n) - (1/2+\Gl_n) \check{f}(n)|^2 \lesssim \| f \|_{\Hcal}^2 + \| g \|_{\Hcal^*}^2
$$
and
$$
\| \Gy \|_{\Hcal^*}^2 \lesssim \| f \|_{\Hcal}^2 + \| g \|_{\Hcal^*}^2.
$$
This proves (ii).

If $\Gl(\Ge_c) = \Gl_n$ for some $n \neq 0$, then we have
$$
\frac{1}{|( \Ge_c -1 ) (\Gl_n-\Gl(\Ge_c)) + i\Gd (\Gl_n-\frac{1}{2})|} \gtrsim \Gd^{-1}
$$
Therefore we have
$$
\| \Gvf^\prime \|_{\Hcal^*} \ge |\hat{\Gvf}(n) | \gtrsim \Gd^{-1} |a_n| .
$$
This completes the proof.
\end{proof}

The following is the main theorem of this paper in three dimensions.

\begin{theorem} \label{mainthm3D}
Suppose $d=3$ and assume
\beq
s \Gd^{-1} \le c
\eeq
for sufficiently small $c$. Let $u_\Gd$ be the solution to \eqnref{eq:Helmholtz2}.
\begin{itemize}
\item[(i)] If $\Gl(\Ge_c/\Ge_m) \neq \Gl_n$ for any $n$, then there is $C$ independent of $\Gd$ (may depend on $\Ge_c/\Ge_m$) such that
\beq
\| \nabla u_\Gd \|_{L^2(\GO)} \le C.
\eeq
\item[(iii)] If $\Gl(\Ge_c/\Ge_m) = \Gl_n$ for some $n \neq 0$, let $z$ be such that $a \cdot \nabla \Scal[\Gvf_n](z) \neq 0$. Then
\beq
\left\Vert \nabla u_\delta \right\Vert_{L^2(\Omega)} \approx \delta ^{-1}
\eeq
as $\delta \to +0$.
\end{itemize}
\end{theorem}

\begin{proof}
We still assume $\Ge_m=1$. Since $A_\Gd^s = A_\Gd + T_\Gd^s= A_\Gd (I + (A_\Gd)^{-1}T_\Gd^s)$, it follows from \eqnref{short} that
$$
\Phi_\Gd = (I + (A_\Gd)^{-1}T_\Gd^s)^{-1} (A_\Gd)^{-1} [\BF],
$$
where
$$
\Phi_\Gd = \begin{bmatrix}
    \Gvf_\Gd \\
    \Gy_\Gd
  \end{bmatrix}
\quad \mbox{and} \quad
\BF  =
  \begin{bmatrix}
    F_z \\ \p_\nu F_z
  \end{bmatrix}.
$$
We see from \eqnref{TGd} and Corollary \ref{cor3D} (i) that
$$
\Vert  (A_\Gd)^{-1}T_\Gd^s \Vert_{\Lcal(X)} \lesssim \Gd^{-1} s.
$$
So, if $s \Gd^{-1}$ is sufficiently small, then we have
\beq\label{Phiest}
\| \Phi_\Gd - (A_\Gd)^{-1} [\BF] \|_X \lesssim \Gd^{-1} s \| (A_\Gd)^{-1} [\BF] \|_X .
\eeq

If $\Gl(\Ge_c) \neq \Gl_n$ for any $n$, then we infer from Corollary \ref{cor3D} (ii) that
$$
\| \Phi_\Gd \|_{X} \le C \| \BF \|_Y.
$$
So, we have from \eqnref{eq:21} and \eqnref{eq:energy_estimate_3d}
$$
\| \nabla u_\Gd \|_{L^2(\GO)} = \| \nabla \Scal^{k_c}[\Gvf_\Gd] \|_{L^2(\GO)} \lesssim \| \Gvf_\Gd \|_{\Hcal^*} \le C
$$
regardless of $\Gd$.

Suppose that $\Gl(\Ge_c) = \Gl_n$ for some $n \neq 0$. Let $(A_\Gd)^{-1} [\BF] = (\Gvf_1, \Gy_1)^T$. Then Corollary \ref{cor3D} (iii) shows that
$$
\| \Gvf_1^\prime \|_{\Hcal^*} \gtrsim |a_n| \Gd^{-1},
$$
where
\beq\label{an}
a_n= (\widehat{\p_\nu F_z})(n) - (1/2+\Gl_n) \check{F_z}(n).
\eeq
It then follows from \eqnref{Phiest} that
$$
\| \Gvf_\Gd^\prime \|_{\Hcal^*} \gtrsim \| \Gvf_1^\prime \|_{\Hcal^*} - \Gd^{-1} s \| (A_\Gd)^{-1} [\BF] \|_X \gtrsim |a_n| \Gd^{-1}
$$
if $a_n \neq 0$ for sufficiently small $\Gd$. Thus we obtain from \eqnref{eq:energy_estimate_3d} that
\beq
\| \nabla u_\Gd \|_{L^2(\GO)} = \| \nabla \Scal^{k_c}[\Gvf_\Gd] \|_{L^2(\GO)} \gtrsim |a_n| \Gd^{-1} - s |\hat{\Gvf_\Gd}(0)|.
\eeq

We now show that $|\hat{\Gvf_\Gd}(0)|$ is bounded, and $a_n \neq 0$ for generic $z$'s. For that purpose we write $A_\Gd^\Go$ as $A_\Gd^\Go = ( I + T_\Gd^s (A_\Gd^0)^{-1} ) A_\Gd^0$ so that \eqnref{short} takes the form
\beq\label{short2}
A_\Gd^0[\Phi_\Gd] = ( I + T_\Gd^s (A_\Gd^0)^{-1} )^{-1}[\BF]
\eeq
Let $( I + T_\Gd^s (A_\Gd^0)^{-1} )^{-1}[\BF] = (f,g)^T$. Then since $\| T_\Gd^s (A_\Gd^0)^{-1} \|_{\Lcal(Y)} \lesssim \Gd^{-1} s$, we have $\| f \|_{\Hcal} + \| g \|_{\Hcal^*}$ is bounded.
Since $\Gl_0=1/2$, we have according to \eqnref{solGvf3D}
$$
|\hat{\Gvf_\Gd} (0)|= \left| \frac{\hat{g}(0) - \check{f}(0)}{( \Ge_c -1 ) (\frac{1}{2}-\Gl(\Ge_c))} \right| \le C.
$$

Recall that $F_z(x) := - a \cdot \nabla_x \GG^{\Go}(x - z)$. According to \eqnref{an} we have
\begin{align*}
a_n &= \la \p_\nu F_z, \Gvf_n \ra_{\Hcal^*} - (1/2+\Gl_n) \la F_z, \Scal[\Gvf_n] \ra_{\Hcal} \\
&= - \la \p_\nu F_z, \Scal[\Gvf_n] \ra + (1/2+\Gl_n) \la F_z, \Gvf_n \ra \\
&= \Go^2 \int_{\GO} F_z \, \Scal[\Gvf_n] \, dx - \la F_z, \p_{\nu}\Scal[\Gvf_n] \big|_{-} \ra + (1/2+\Gl_n) \la F_z, \Gvf_n \ra \\
&= \Go^2 \int_{\GO} F_z \, \Scal[\Gvf_n] \, dx  + \la F_z, \Gvf_n \ra.
\end{align*}
Since $F_z(x) = a \cdot \nabla_z \GG^{\Go}(x - z)$, we have
$$
\la F_z, \Gvf_n \ra = a \cdot \nabla \Scal^{\Go} [\Gvf_n](z).
$$
By \eqnref{eq:10} we have
$$
\nabla \Scal^{\Go} [\Gvf_n](z) = \nabla \Scal [\Gvf_n](z) + O(\Go^2),
$$
and hence
$$
a_n = a \cdot \nabla \Scal [\Gvf_n](z) + O(\Go^2)
$$
Note that $a \cdot \nabla \Scal [\Gvf_n](z)$ is a harmonic function in $z \in \Rbb^3 \setminus \ol{\GO}$. So it cannot be zero for $z$ in an open set. We choose $z$ so that $a \cdot \nabla \Scal [\Gvf_n](z) \neq 0$, and then $a_n \neq 0$ if $\Go$ is sufficiently small.
Thus we have
\beq
\| \nabla u_\Gd \|_{L^2(\GO)} \gtrsim \Gd^{-1} .
\eeq
This completes the proof.
\end{proof}

\subsection{Two dimensions}\label{sec:2d}

In two dimensions we decompose $A_{\Gd}^s$ in \eqref{eq:A_delta^omega} into three parts: $A_\Gd^s = A_\Gd + B^s + T_\Gd^s$ where $A_\Gd$ is defined by \eqnref{A_delta^0} and
\beq\label{Bdef}
 B^s =
  \begin{bmatrix}
    \Gj^{k_c} \left\la \cdot, 1 \right\ra & - \Gj \left\la \cdot, 1 \right\ra \\
    0 & 0
  \end{bmatrix} .
\eeq
Here, $\Gj^{k_c}$ is defined by
\beq
  \Gj^{k_c} = \left( 1 / 2 \pi \right) \left( \ln {k_c} + \Gg - \ln {2} \right) - i / 4,
\eeq
and $\Gj$ is defined by \eqref{eq:tau}. We emphasize that
\beq
  \vert \Gj^{k_c} \vert \sim - \ln {\Go}, \quad
  \vert \Gj \vert \sim - \ln {\Go}. \label{eq:5.2}
\eeq
We have from \eqnref{Stwo}, \eqnref{Qtwo} and Lemma \ref{regularity}
\beq\label{TGd2}
  \left\Vert T_\Gd^s \right\Vert_{\Lcal(X, Y)} \lesssim |s^2 \ln s|.
\eeq

Unlike the three dimensional case, $A_\Gd: X \to Y$ may not be invertible since $\Scal: \Hcal^* \to \Hcal$ is not invertible in general. Instead we prove that $A_\Gd + B^s: X \to Y$ is invertible. In fact, we obtain the following lemma.

\begin{lemma}\label{thm2D}
The operator $A_\Gd + B^s: X \to Y$ is invertible. For $(f,g)^T \in Y$, the solution $(\Gvf,\psi)^T$ to the equation
$$
(A_\Gd + B^s) \begin{bmatrix} \Gvf \\ \Gy \end{bmatrix} = \begin{bmatrix} f \\ g \end{bmatrix}
$$
is given by
\beq\label{solGvf2D}
\Gvf=\Gvf^\prime + \hat{\Gvf}(0) \Gvf_0 = \Gvf^\prime + \frac{\check{f}(0) \tilde{\Scal}[\Gvf_0] - \hat{g}(0) \left( \Scal[\Gvf_0] + \Gj \la \Gvf_0, 1 \ra \right)}{\Scal[\Gvf_0] + \Gj^{k_c} \la \Gvf_0, 1 \ra} \Gvf_0
\eeq
and
\beq\label{solGy2D}
\Gy = \Gvf^\prime - \Scal^{-1}[f^\prime] - \hat{g}(0) \Gvf_0,
\eeq
where
\beq\label{Gvfprime}
\Gvf^\prime = \sum_{n=1}^\infty \frac{\hat{g}(n) - (\frac{1}{2}+\Gl_n) \check{f}(n)}{( \Ge_c -1 ) (\Gl_n-\Gl(\Ge_c)) + i\Gd (\Gl_n-\frac{1}{2})} \Gvf_n.
\eeq
\end{lemma}

Before proving Theorem \ref{thm2D}, we emphasize that $\Scal[\Gvf_0]$ is constant ($=c_0$). If $c_0 \neq 0$, then $\tilde{\Scal}[\Gvf_0] =  \Scal[\Gvf_0]=c_0$, and
$$
\la \Gvf_0, 1 \ra = c_0^{-1} \la \Gvf_0, \Scal[\Gvf_0] \ra = c_0^{-1}.
$$
So we have
$$
\hat{\Gvf}(0)= \frac{c_0 \check{f}(0) - (c_0+ c_0^{-1}\Gj) \hat{g}(0) }{c_0 + c_0^{-1} \Gj^{k_c}} .
$$
If $c_0 =0$, then $\tilde{\Scal}[\Gvf_0] = 1$, and $\la \Gvf_0, 1 \ra = 1$.
So we have
$$
\hat{\Gvf}(0)= \frac{\check{f}(0) - \Gj \hat{g}(0)}{\Gj^{k_c}}.
$$

\noindent{\sl Proof of Lemma \ref{thm2D}}.
Let
$$
f=f^\prime + \check{f}(0) \tilde{\Scal}[\Gvf_0] \quad g=g^\prime + \hat{g}(0) \Gvf_0
$$
be orthogonal decompositions in $\Hcal$ and $\Hcal^*$, so that $f^\prime \in \Hcal_0$ and $g^\prime \in \Hcal_0^*$.

Since $\Scal: \Hcal_0^* \to \Hcal_0$ is invertible, one can see as in Lemma \ref{thm3D} that
the solution to
$$
A_\Gd \begin{bmatrix} \Gvf^\prime \\ \Gy^\prime \end{bmatrix} = \begin{bmatrix} f^\prime \\ g^\prime \end{bmatrix}
$$
is given by \eqnref{Gvfprime} and
$$
\Gy^\prime = \Gvf^\prime - \Scal^{-1}[f^\prime].
$$

Since $(-1/2I + \Kcal^*)[\Gvf_0]=0$ and $\Gvf^\prime, \Gy^\prime \in \Hcal_0^*$, we can see that
\begin{align*}
(A_\Gd + B^s)
\begin{bmatrix} \Gvf^\prime + c \Gvf_0 \\ \Gvf^\prime - \Scal^{-1}[f^\prime] + d \Gvf_0 \end{bmatrix}
&= \begin{bmatrix} f^\prime \\ g^\prime \end{bmatrix}
+ (A_\Gd + B^s) \begin{bmatrix} c \Gvf_0 \\ d \Gvf_0 \end{bmatrix} \\
&= \begin{bmatrix} f^\prime \\ g^\prime \end{bmatrix}
+ \begin{bmatrix} c \left( \Scal[\Gvf_0] + \Gj^{k_c} \la \Gvf_0,1 \ra \right) - d \left( \Scal[\Gvf_0] + \Gj \la \Gvf_0,1 \ra \right) \\
- d \Gvf_0 \end{bmatrix}.
\end{align*}
So we solve
$$
c \left( \Scal[\Gvf_0] + \Gj^{k_c} \la \Gvf_0,1 \ra \right) - d \left( \Scal[\Gvf_0] + \Gj \la \Gvf_0,1 \ra \right) = \check{f}(0)\tilde{\Scal}[\Gvf_0], \quad -d = \hat{g}(0)
$$
for $c,d$ to have \eqnref{solGvf2D} and \eqnref{solGy2D}. This completes the proof.
\qed

We can obtain from Lemma \ref{thm2D} results similar to those in Corollary \ref{cor3D} for two dimensions. We then obtain the following theorem for two dimensions.

\begin{theorem} \label{mainthm2D}
Suppose $d=2$ and assume
\beq
s^2 |\ln s| \Gd^{-1} \le c
\eeq
for sufficiently small $c$. Let $u_\Gd$ be the solution to \eqnref{eq:Helmholtz2}.
\begin{itemize}
\item[(i)] If $\Gl(\Ge_c/\Ge_m) \neq \Gl_n$ for any $n$, then there is $C$ independent of $\Gd$ (may depend on $\Ge_c/\Ge_m$) such that
\beq
\| \nabla u_\Gd \|_{L^2(\GO)} \le C.
\eeq
\item[(iii)] If $\Gl(\Ge_c/\Ge_m) = \Gl_n$ for some $n \neq 0$, let $z$ be such that $a \cdot \nabla \Scal[\Gvf_n](z) \neq 0$. Then
\beq
\left\Vert \nabla u_\delta \right\Vert_{L^2(\Omega)} \approx \delta ^{-1} \quad \mbox{as } \delta \to +0.
\eeq
\end{itemize}
\end{theorem}

\begin{proof}
We write
\beq
A_\Gd^s = A_\Gd + B^s + T_\Gd^s = (A_\Gd + B^s)\left( I + (A_\Gd + B^s)^{-1}T_\Gd^s \right),
\eeq
and follow the same lines of the proof for Theorem \ref{mainthm3D}. One thing we need to check is that $|\hat{\Gvf_\Gd}(0)|$ is bounded. To do that it suffices to show that $|\hat{\Gvf}(0)|$ is bounded where $\Gvf$ is the solution expressed in \eqnref{solGvf2D}. Note that
$$
|\hat{\Gvf}(0)| = \left| \frac{\check{f}(0) \tilde{\Scal}[\Gvf_0] - \hat{g}(0) \left( \Scal[\Gvf_0] + \Gj \la \Gvf_0, 1 \ra \right)}{\Scal[\Gvf_0] + \Gj^{k_c} \la \Gvf_0, 1 \ra} \right| \lesssim \frac{|\Gj^{k_c}|}{|\Gj|} \lesssim 1.
$$
This completes the proof.
\end{proof}

\section*{Acknowledgement}

We thank Habib Ammari for sending us the paper \cite{AMRZ15}. The work of K. Ando and H. Kang was supported by the Korean Ministry of Education, Sciences and Technology through NRF grant No. 2010-0017532. The work of H. Liu was supported by the FRG and start-up grants of Hong Kong Baptist University, and the NSF grant of China, No.\,11371115.



\end{document}